\documentclass[reprint,aps,pre]{revtex4-1}
\usepackage{mathrsfs}
\usepackage{algorithmic}
\usepackage{graphicx}
\usepackage{textcomp}
\usepackage{wrapfig}
\def\xx{\ensuremath\mathbf{x}}
\def\XX{\ensuremath\mathbf{X}}

\def\sigsig{\ensuremath\mathbf{\Sigma}}
\def\VV{\ensuremath\mathbf{V}}

\def\sss{\ensuremath\mathbf{s}}

\def\pp{\ensuremath\mathbf{\Phi}}
\def\psps{\ensuremath\mathbf{\Psi}}
\def\bemat{\ensuremath\begin{bmatrix}}
\def\eemat{\ensuremath\end{bmatrix}}

\def\aa{\ensuremath\mathbf{a}}

\usepackage{amsmath, amssymb, amscd, amsthm, amsfonts}
\usepackage{graphicx}
\usepackage{hyperref}
\DeclareMathOperator*{\argmax}{arg\,max}
\DeclareMathOperator*{\argmin}{arg\,min}
\newcommand{\rr}{\mathbb{R}}

\begin{document}
\title{Data-driven sensor placement with shallow decoder networks}
\author{
Jan~Williams$^*$,  Olivia~Zahn$^\dag$, J.~Nathan~Kutz$^\ddag$\\[.1in]
{\em \small
$^*$ {Departments of Physics and Mathematics, Carleton College, Northfield, MN 55057 USA}\\
$^\dag$ Department of Physics, University of Washington, Seattle, WA, 98195 USA\\
$^\ddag$ Department of Applied Mathematics, University of Washington, Seattle, WA, 98195 USA }
}

\begin{abstract}
Sensor placement is an important and ubiquitous problem across the engineering and physical sciences for tasks such as reconstruction, forecasting and control.
Surprisingly, there are few principled mathematical techniques developed to date for optimizing sensor locations, with the leading sensor placement algorithms often based upon the discovery of linear, low-rank sub-spaces and the QR algorithm.
QR is a computationally efficient greedy search algorithm which selects sensor locations from candidate positions with maximal variance exhibited in a training data set.
More recently, neural networks, specifically {\em shallow decoder networks} (SDNs), have been shown to be very successful in mapping sensor measurements to the original high-dimensional state space.
SDNs outperform linear subspace representations in  reconstruction accuracy, noise tolerance, and robustness to sensor locations.
However, SDNs lack principled mathematical techniques for determining sensor placement.
In this work, we develop two algorithms for optimizing sensor locations for use with SDNs:  one which is a linear selection algorithm based upon QR (Q-SDN), and one which is a nonlinear selection algorithm based upon neural network pruning (P-SDN).
Such sensor placement algorithms promise to enhance the already impressive reconstruction capabilities of SDNs.
We demonstrate our sensor selection algorithms on two example data sets from fluid dynamics.  Moreover, we provide a detailed comparison between our linear (Q-SDN) and nonlinear (P-SDN) algorithms with traditional linear embedding techniques (proper orthogonal decomposition) and QR greedy selection.
We show that QR selection with SDNs enhances performance.  QR even out-performs our nonlinear selection method that uses magnitude-based pruning.
Thus, the combination of a greedy linear selection (QR) with nonlinear encoding (SDN) provides a synergistic combination.
\end{abstract}


\maketitle

\begin{figure*}[!t]
\centerline{\includegraphics[width=\textwidth]{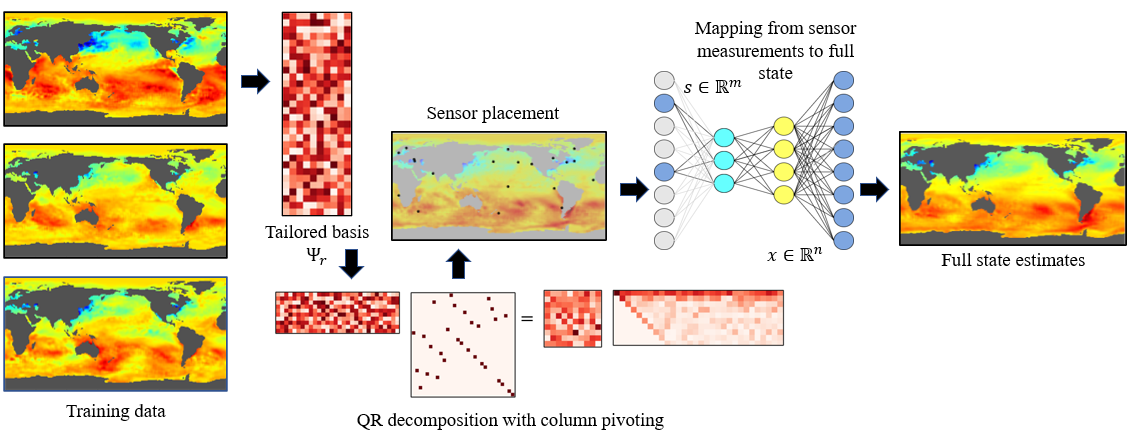}}
\caption{Flow diagram of our data-driven method for high-dimensional state estimation, Q-SDN.  We first obtain a tailored basis $\mathbf{\Psi} _r$ for our training data via a singular value decomposition.  We then perform a QR decomposition with column pivoting to obtain pivots corresponding to our chosen sensor locations.  Finally, we train a shallow neural network to learn an end-to-end, nonlinear mapping from data ``collected'' at these sensor locations ($\sss \in \mathbb{R}^m$) to the full-state ($\xx \in \mathbb{R}^n$).}
\label{fig:flow_diagram}
\end{figure*}

\section{Introduction}
Throughout engineering and the physical sciences, we rely on discrete measurements to infer the current and future behavior of spatio-temporal systems.  There are often significant limitations in the number of measurements possible when studying a particular system.  For instance, additional sensors might be prohibitively expensive or otherwise infeasible to place in certain spatial locations.  Determining ideal sensor locations can greatly enhance the ability to reconstruct, forecast or control a system under consideration.  Despite the ubiquity of this problem, there remain relatively few principled mathematical methods for optimizing sensor placement.  The QR algorithm is one version of a greedy search that is computationally tractable, easy to implement, and demonstrates robust performance.  Reconstructions with limited (sparse) sensors typically rely on computing low-rank embeddings of the spatio-temporal dynamics.  Dimensionality reduction is commonly achieved with the {\em singular value decomposition} (SVD), which is otherwise known as the {\em proper orthogonal decomposition} (POD) and closely related to {\em principal component analysis} (PCA)~\cite{kutz2013data,brunton2019data}.  The dominant correlated modes of the SVD provide a linear map between measurements and the full spatio-temporal state space~\cite{manohar-sensors}.  However, sensor placement and reconstruction using QR with SVD is a linear method which fails to leverage recent developments in the descriptive power of nonlinear techniques and deep learning.  In this work, our goal is to identify sensor locations with two algorithms that use neural-network based reconstructions: (i) we identify sensor locations using linear embeddings and the QR algorithm, but with a shallow decoder for reconstruction, and (ii) we identify the nonlinear selection of sensor locations by pruning (systematic removal of weights/nodes) an autoencoder architecture to a shallow decoder network.  The two methods are compared against standard linear sensor placement and reconstruction methods (QR and SVD), showing that by exploiting the nonlinearity of deep learning, the performance of the sensor selection problem can be made more accurate and robust to noise.  Additionally, in the low sensor limit, the difference in reconstruction performance between the shallow decoder with optimal sensors and traditional techniques is greater. 
%

For moderate-sized search spaces, the sensor placement problem has well-known model-based solutions using optimal experiment design~\cite{boyd2004convex,joshi2008sensor}, and information theoretic and Bayesian criteria~\cite{caselton1984optimal,krause2008near,lindley1956measure,sebastiani2000maximum,paninski2005asymptotic}.  As an alternative, greedy algorithms can provide fast, near-optimal candidate locations.  Indeed, many of the current linear methods for reconstructing full-state spaces from sparse measurements are based on the mathematical framework of gappy POD \cite{gappy-pod} which has been demonstrated to be successful in applications ranging from ocean modeling \cite{manohar-sensors, ocean-modeling} to aerodynamics \cite{aerodynamics}.  Gappy POD and its derivative techniques function by estimating inner products from discrete measurements, a process described in more detail in Section \ref{section:placement}.  While gappy POD was initially used with random measurements, subsequent work has proposed principled methodologies for optimizing measurement selection to estimate these inner products; the empirical interpolation method (EIM) \cite{EIMS}, discrete empirical interpolation method (DEIM) \cite{DEIMS}, and a variant of DEIM using a QR decomposition (Q-DEIM) \cite{Q-DEIMS} offer protocols for determining ideal measurements.  Manohar et al~ \cite{manohar-sensors} and Peherstorfer et al~\cite{peherstorfer2020stability} demonstrate the efficacy of these principled  selection criteria for use with linear, interpolative reconstructions.  Such methods can also easily incorporate cost constraints~\cite{clark2018greedy,clark2020sensor}, multi-fidelity sensors~\cite{clark2020multi}, multi-scale physics~\cite{manohar2019optimized}, and control requirements~\cite{manohar2018optimal}. Importantly, if a sufficiently large number of measurements is available, then random sampling and reconstruction with SVD-based modes works just as well as the principled selection methods.  In this case, sensor selection is simply not an issue.  However, in the limit of a small number of sensors, the reconstruction error using random sampling increases sharply.  The minimum limit on the number of sensors is problem and application specific and depends upon the underlying low-rank spatio-temporal dynamics.  It is in this limit that greedy algorithms are critically enabling as a brute-force combinatorial search for optimal placement is computationally infeasible.  

Machine learning and AI are changing the ways in which we can perform reconstructive tasks~\cite{goodfellow2016deep,brunton2019data}, including with optimal sensors~\cite{semaan2017optimal}.  More recently in the context of sensors and reconstruction, Erichson et al \cite{erichson-decoder} demonstrated the ability of {\em shallow neural networks} (SDNs) to infer high-dimensional state spaces from few sensor measurements.  SDNs use a set of training data to learn a nonlinear mapping from measurements to high-dimensional state spaces.  The trained networks can then be used to reconstruct the full state space from a dramatically limited subsampling.  SDNs enjoy a variety of advantages when compared to traditional, linear methods.  Namely, reconstructions using SDNs have been demonstrated to exceed the performance of linear reconstructions in accuracy, noise tolerance, and robustness to sensor locations \cite{erichson-decoder}.  Here, we demonstrate that although SDNs are less sensitive to sensor location than linear methods, the reconstructive performance of SDNs can be further improved by selecting measurement locations in a principled manner.  We develop two sensor selection methods (one linear and one nonlinear) and combine these choices of sensors with the SDN.  
Importantly, the architecture advocated provides performance gains for a minimal number of sensors.  Such scenarios are commonly considered when the sensors, or placing them, is prohibitively expensive (e.g.  ocean or atmospheric monitoring).  

The first selection algorithm we consider for use with SDNs follows that of \cite{manohar-sensors, Q-DEIMS} and uses an extension of Q-DEIM to identify sensor locations by selecting locations with the highest variance in POD modes obtained through a singular value decomposition of a matrix of training data.  In essence, we combine a linear sensor selection (QR) method with a nonlinear shallow decoder reconstruction (SDN) to produce our Q-SDN algorithm, demonstrated in Fig.  \ref{fig:flow_diagram}.  This selection methodology is attractive for a number of reasons.  Primarily it is easy and computationally inexpensive to compute and has been demonstrated to identify physically meaningful patterns in data \cite{manohar-sensors}.  Empirically, we also find it yields superior performance in comparison to our nonlinear selection method.  In summary, Q-SDN is a linear selection with nonlinear reconstruction.



In addition to Q-SDN, we also implement neural network pruning as a nonlinear selection method.  Neural network pruning is a biologically inspired method for sparsifying network architecture by systematically removing weights or nodes.  There are many methods for choosing which weights/nodes to prune.  Generally, the methods either evaluate the effect of a weight/node on the performance (i.e.  loss) or prediction \cite{NIPS1989_6c9882bb, hassibi1993optimal}.  For example, simple magnitude-based weight pruning removes weights whose absolute values are closest to zero.  Along with other techniques (i.e.  regularization during training), pruning has been shown to be an effective method for sparsifying a DNN without compromising performance \cite{NIPS1989_6c9882bb, hassibi1993optimal, louizos2018learning, louizos2017bayesian, kuzmin2019taxonomy}.  Pruning, like dropout \cite{JMLR:v15:srivastava14a}, can also help prevent overfitting and improve the generalization error.  When applied to the input layer of a network, pruning can identify important input features by removing the inputs that contribute minimally to the prediction \cite{olivia}.  Here we use pruning as a method for discovering a set of sparse sensors by pruning only the input layer of an autoencoder network.  The method begins by training an autoencoder to map full state sensor measurements back to the original high-dimensional state space.  The first layer of the network is the encoder network which maps the full state sensor measurements to a low-dimensional representation.  Following the method in \cite{olivia}, we prune this layer iteratively.  Specifically, once the training error plateaus, we prune a small subset of input nodes with corresponding weights that have the smallest root mean square.  The pruned network is then retrained  and this process repeats until only a few sensor measurements remain.  This results in a pruned version of the encoder, while the decoder remains fully connected.  This pruned SDN (P-SDN) is highly intuitive in its selection as those locations which have the smallest impact on prediction are removed during training.  The P-SDN simultaneously learns sensor locations and reconstructions, which is quite different than Q-SDN which learns these in serial.

In Section \ref{section:placement}, we introduce a linear method for reconstructing high-dimensional states and our linear placement algorithm.  In Section \ref{section:decoders}, we formally present how SDNs are used to perform reconstructions, both Q-SDN and P-SDN, and we discuss how identifying sensor placement by neural network pruning is a natural extension of SDNs.  Section \ref{section-experiments} compares the performance of our placement algorithms in two datasets, one synthetic and one natural, and finally, in Section \ref{section:discussion} we discuss our results.  



\section{Low-Rank Reconstructions and QR Sensor Placement}
\label{section:placement}

Traditional data-driven techniques for reconstructing states given a limited number of sensors often rely on estimating the state through low-rank approximations.  Specifically, we use POD modes, which are the left eigenvectors of the singular value decomposition of a matrix.  The POD modes give the dominant correlated spatial dynamics of the data.  Thus the columns consist of orthonormal, representative samples of the state.  Sensor data is then used to estimate the coefficients of these POD modes for new samples.  In this section, we outline this process and introduce the QR decomposition with column pivoting as a means to determine a greedy approximation for optimal sensor placement.

\subsection{Low-Rank Reconstructions}
\label{section:POD}
Let $\xx \in \rr^m$ represent a sample state of a physical system.  We seek to estimate the full state $\xx$ given sensor measurements $\sss \in \rr ^n$ with $n \ll m.$ We restrict our consideration to the case that sensors perform discrete measurements (e.g.  measuring some quantity at a particular location).  In this case, we have $s _j = x _i$ for some $1 \leq j \leq n$ and $1 \leq i \leq m.$ Furthermore, we can write
\begin{equation}
    \pp \xx = \sss
    \label{eq:sensing}
\end{equation}
for some $n \times m$ matrix $\pp$ constructed by taking $n$ rows of the $m \times m$ identity matrix.\\
\indent In general, Eq.  (\ref{eq:sensing}) represents an underdetermined linear system and, given $\sss$ and $\pp$, reconstructing $\xx$ represents an ill-posed problem. However, we can leverage the fact that many physical systems have low dimensional approximations to determine an approximate solution for Eq.  (\ref{eq:sensing}).  To do so, we require a set of training data, $\{ \xx_i \in \rr ^m \}_{i=1}^N$, consisting of $N$ representative samples of the physical state.  From this training data, we construct the data matrix 
\begin{equation}
    \XX = \bemat | &  & | \\ \xx_1 & \dots & \xx_N\\  | &  & |\eemat.
\end{equation}
Then, the singular value decomposition of $\XX$ yields 
\begin{equation}
    \XX = \mathbf{\Psi \Sigma V}^T
    \label{eq:datamatrix}
\end{equation}
where $\psps \in \rr ^{m \times m}$ and $\VV \in \rr ^{n \times n}$ are orthogonal matrices and $\sigsig \in \rr ^{m \times n}$ is a diagonal matrix with non-increasing entries.  The optimal rank $r$ approximation of $\XX$ is given by 
\begin{equation}
    \XX \approx \hat{\XX} = \psps _r \sigsig _r \VV ^T _r = \sum _{i = 1} ^r \sigma _i \mathbf{\psi}_i v_i^T
\end{equation}
where $\sigma _i$ is the $(i,i)$ entry in $\sigsig$, $\mathbf{\psi} _i$ the $i$-th column in $\psps$, and $v _i$ the $i$-th column in $\VV$.  From this, we see that the $i$-th column of $\XX$ is approximated by a linear combination, given by $\aa_i \in \rr ^r$, of the columns of $\psps _r$.  That is,
\begin{equation}
    \xx _i \approx \hat{\xx}_i := \psps _r \aa_i.
\end{equation}
Since we assume the columns of $\XX$ to be representative of the same physical system as $\xx,$ there also exists $\aa \in \rr ^r$ such that 
\begin{equation}
    \xx \approx \hat{\xx} = \psps_r \aa.
    \label{eq:approx_x}
\end{equation}
Combining this with Eq.  (\ref{eq:sensing}) gives, 
\begin{equation}
    \sss \approx \pp \psps_r \aa.
    \label{eq:approx_a}
\end{equation}
Eq.  (\ref{eq:approx_a}) serves as the basis of POD reconstruction techniques.  By selecting the number of sensor measurements, $n$, to be greater than or equal to the number of modes used to estimate the data, $r$, we can determine $\aa$ by
\begin{equation}
    \aa = \begin{cases} 
      (\pp \psps_r)^{-1} \sss, \; n = r \\
      (\pp \psps_r)^{\dagger}\sss, \; n > r\\
   \end{cases},
   \label{eq:a_cases}
\end{equation}
where $\dagger$ denotes the Moore-Penrose pseudoinverse.  Subsequently, applying Eq.  (\ref{eq:approx_x}) yields an estimate of the state $\xx.$\\

\subsection{QR Decomposition with Column Pivoting}
The previous section describes a method for estimating a full state $\xx$ given some sub-sampling $\sss,$ but does not consider the influence of which entries of $\xx$ are selected.  That is, it does not consider the structure of $\pp,$ or, returning to the physical interpretation of $\xx$, ``where'' in the state measurements are taken.  In many applications, however, the placement of sensors can drastically impact one's ability to accurately infer the full state; for instance, taking all measurements in only one region of a field might be less effective than sampling throughout the field.  More formally, we seek a matrix $\pp$ that leaves the inversion from Eq.  (\ref{eq:a_cases}) well-conditioned.  The following argument \cite{manohar-sensors} justifies the use of a QR matrix decomposition as a greedy approximation for an optimal $\pp$ when the number of modes used in the reconstruction is equal to the number of measurements, $n = r$.  The argument can be extended to the oversampled case, $n > r$, but the remainder of our considerations, including the computational experiments, assume $n = r.$

\begin{figure}[t]
    \centering
    \includegraphics[width=\columnwidth]{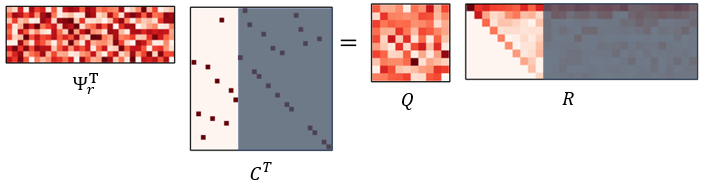}
    \caption{The decomposition of a tailored basis $\psps _r ^T$ into a column permutation matrix $\mathbf{C^T}$, orthogonal matrix $\mathbf{Q}$ and a upper triangular matrix $\mathbf{R}$ with non-increasing entries on the diagonal.  The non-shaded blocks represent $\mathbf{C}_r^T$ (left) and $R_r$ (right).}
    \label{fig:qr_placement}
\end{figure}

\indent Recall that the rows of $\pp$ will have exactly one non-zero entry, equal to one.  Let $\gamma_i \in \mathbb{N}$ with $1 \leq \gamma _i \leq m$ denote the the non-zero entry of the $i$-th row of $\pp$.  Now, let $\Gamma = \{ \gamma _i \; : 1 \leq i \leq r\; \}$, and let $\pp _\Gamma$ be the matrix uniquely determined by $\Gamma$.  An indirect bound on the condition number of $(\pp _\Gamma \psps _r)$ is found by 
\begin{align}
    \pp _\Gamma  .  = \argmax_{\pp _\Gamma, |\Gamma|=r} |\det \pp _\Gamma \psps _r| = \argmax_{\pp _\Gamma, |\Gamma|=r} \prod _i|\lambda _i (\pp _\Gamma \psps _r)| \notag \\ = 
    \argmax_{\pp_\Gamma, |\Gamma|=r} \prod _i\sigma _i (\pp _\Gamma \psps _r).
    \label{eq:qr_optimization}
\end{align}
Unfortunately, direct optimization is computationally intractable and requires an NP-hard, combinatorial search.  However, a greedy approximation is given by the QR decomposition with column pivoting of $\psps _r ^T$.  The decomposition finds a column permutation matrix $\mathbf{C}^T$, an orthogonal matrix $\mathbf{Q}$, and an upper triangular matrix $\mathbf{R}$ such that $\psps_r ^T \mathbf{C}^T = \mathbf{QR}$.  The diagonal entries of $\mathbf{R}$ are determined by selecting the pivot column with maximum $\ell_2$ norm.  The orthogonal projection of the pivot column is then subtracted from all other columns and the process repeats, ensuring that each diagonal entry of $\mathbf{R}$ will have maximal $\ell_2$ norm and thus be non-increasing.  Now, consider the matrix constructed by truncating $\mathbf{C}^T$ to include only the first $r$ columns, call it $\mathbf{C}^T_r$, and take the product of this submatrix with $\mathbf{\psps}_r^T.$ Now, we have $\psps_r^T \mathbf{C}_r^T = \mathbf{Q}\mathbf{R}_r,$ where $\mathbf{R}_r$ denotes the matrix constructed by taking the first $r$ columns of $\mathbf{R}$.  Note that we have a greedy approximation for the maximum volume of this matrix product since, by construction, each diagonal entry of $\mathbf{R}$ is chosen to be as large as possible, given the preceding entries.  Then, since $\mathbf{R}_r$ is upper triangular and square and $\mathbf{Q}$ is orthogonal,
\begin{equation}
    |\det \psps _r ^T \mathbf{C}_r^T| = \prod _i ^r |r_{ii}|.
\end{equation}
Furthermore, 
\begin{equation}
    \det \psps _r ^T \mathbf{C}_r^T = \det \mathbf{C}_r \psps _r 
\end{equation}
so $\pp _\Gamma = \mathbf{C}_r$ is a greedy approximation of the optimal sensing matrix from Eq.  (\ref{eq:qr_optimization}).  The matrix diagram in Fig.  \ref{fig:qr_placement} pictorially represents this decomposition.\\
\indent Previous work has shown that the greedy sensor placement outlined above significantly outperforms random placement in a variety of reconstructions performed as outlined in Section \ref{section:POD} \cite{manohar-sensors, Q-DEIMS,peherstorfer2020stability} and identifies physically meaningful features, such as the El Nino in sea surface temperature data \cite{manohar-sensors}.  The combination of this measurement selection with SDNs is our proposed Q-SDN algorithm.\\

\section{Shallow Decoders for Reconstruction}
\label{section:decoders}
The previous section provides both a methodology for reconstructing high-dimensional states from some set of measurements as well as an algorithm for identifying near-optimal sensor locations.  However, the work of \cite{erichson-decoder} has demonstrated that shallow neural networks can perform superior reconstructions when compared to linear methods.  A novel contribution of our work is the combination of the linear placement from the previous section with these decoder networks and a nonlinear selection method based on neural network pruning.\\
\subsection{Shallow Decoder Networks (SDNs)}
\indent Let 
\begin{equation}
    \pp \xx = \sss
\end{equation}
represent a sampling of a high-dimensional state.  SDNs are shallow, fully-connected neural networks that have learned a mapping from this subsampling back to the high-dimensional state.  We emphasize the use of a shallow neural network for a few reasons; shallow neural networks tend to be easier to train with limited training data and require less tuning of hyperparameters (citations here).\\
\indent The fully-connected neural networks we use can be denoted as 
\begin{equation}
    \mathcal{F}(\sss; \mathbf{W}) := R(\mathbf{W}^kR(\mathbf{W}^{k-1}\cdots R(\mathbf{W^1 s}))),
\end{equation}
where $\sss$ is the input sensor data, $R$ is a scalar, nonlinear activation function, each $\mathbf{W^i}$ is a weight matrix, and $k$ represents the number of layers in the network.  Let $\mathscr{F}$ be the set of all such networks.  Formally, the desired SDN is $\mathcal{F} \in \mathscr{F}$ such that reconstruction error is minimized.  That is,
\begin{equation}
    \mathcal{F} \in \argmin_{\widetilde{\mathcal{F}} \in  \mathscr{F}} \sum _{i=1}^N ||\xx_i - \widetilde{\mathcal{F}}(\sss_i)||_2
\end{equation}
for a set of training states $\{ \xx _i\}_{i=1}^N$ and corresponding subsamplings of these states, $\{ \sss _i\}_{i=1}^N$.  The architecture described here follows the work of Erichson et al.  \cite{erichson-decoder}.  We use the ADAM optimization algorithm to train our networks, however, we vary the hyperparameters depending on the dataset we consider to achieve the best performance.\\
\indent After training, a new state $\xx$ is approximated from sensor data $\sss$ by evaluating 
\begin{equation}
    \mathcal{F} (\sss) = \hat{\xx}.
\end{equation}
As in Section \ref{section:placement}, a function $\mathcal{F}$ can be learned for measurements obtained through any matrix $\pp$.  Previous work suggests that neural-network based reconstructions are more agnostic to measurement location than linear techniques \cite{erichson-decoder}, however an investigation of optimal sensor placement for neural-network reconstructions has not yet been performed.  We propose combining SDNs with the placement algorithm from Section \ref{section:placement} and also present a novel placement scheme based on neural network pruning in the following subsection.\\
\subsection{Sensor Placement via Pruning}
\begin{figure}[t]
    \centering
    \includegraphics{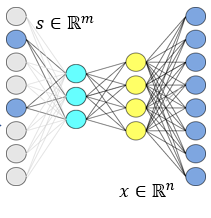}
    \caption{Shallow-decoder network diagram mapping sensor measurements, $\sss \in \rr^m$ to $\xx \in \rr ^n$.  The gray nodes represent elements of the full-state $\xx$ that are not measured, while the blue nodes represent measurements.  The question of optimal sensor placement for reconstruction via a shallow-decoder amounts to determining the optimal active nodes.}
    \label{fig:network_diagram}
\end{figure}

\indent In the context of network-based reconstructions, optimal sensor placement can be viewed as identifying optimal input nodes in a SDN that maps full-states to full-states, depicted in Fig.  \ref{fig:network_diagram}.  In the case that all input nodes are active, the network represents an autoencoder.  Limiting sensor measurements corresponds to removing all connections from an input node in such a network.

\begin{figure}[!t]
    \centering
    \includegraphics[width=\columnwidth]{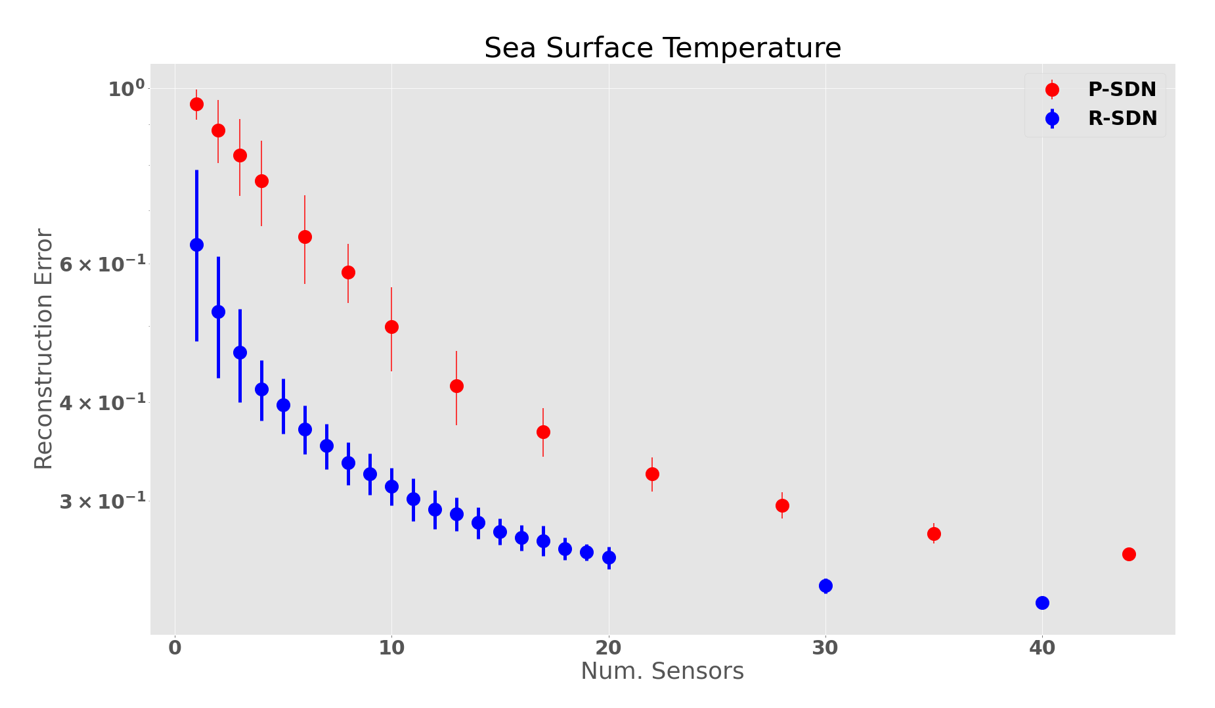}
    \caption{Reconstruction performance of R-SDN and P-SDN for the example dataset of sea surface temperature. Surprisingly, randomly placed sensors outperform those identified by our iterative pruning algorithm.}
    \label{fig:p-sdn-fig}
\end{figure}

\begin{figure}[t]
    \centering
    \includegraphics[width=0.65\columnwidth]{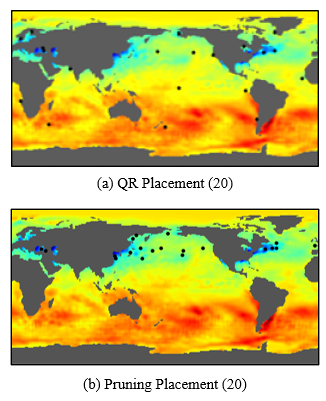}
    \caption{Placement of 20 sensors obtained by a QR decomposition (a) and iterative, magnitude based pruning (b).  QR placement scatters sensors throughout the field, while the placement obtained by pruning clusters sensors in confined geographic regions.}
    \label{fig:pruning}
\end{figure}

\begin{figure*}[!t]
    \centering
    \includegraphics[width=\textwidth]{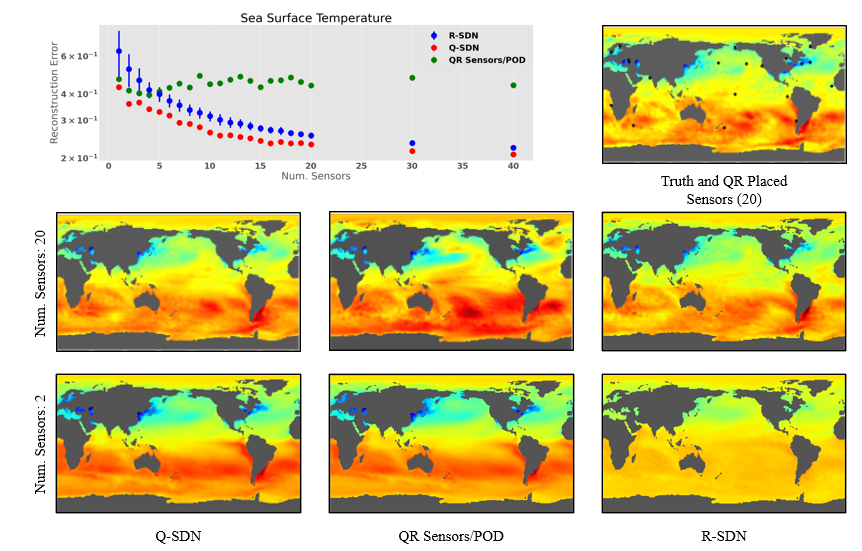}
    \caption{Sea-surface temperature reconstruction errors obtained for varying numbers of sensors and example reconstructions obtained by Q-SDN, QR placed sensors with POD reconstruction, and R-SDN.}
    \label{fig:sst_panel}
\end{figure*}

\indent Neural network pruning is a longstanding technique for sparsifying neural networks by the systematic removal of weights or nodes \cite{NIPS1989_6c9882bb, hassibi1993optimal, baykal2022sensitivity}.  Pruned networks have been shown to retain comparable accuracy to larger models, and as networks have grown in size, interest in pruning as a method for limiting the memory required to store DNNs has increased.  Although pruning can be applied to either weights or nodes, here we consider the removal of nodes.  In this context, we use a form of magnitude-based pruning which identifies nodes for removal by computing the root mean square of the trained weights corresponding to each node.  The nodes associated with the smallest root mean square weights are removed.  Pruning nodes from the input layer of an SDN mapping full-states to full-states has the physical interpretation of removing a sensor measurement.\\
\indent Here, we follow the iterative pruning protocol in \cite{olivia}.  First, we train a fully-connected SDN to convergence.  We define a target sparsity (i.e.  10\%, 20\%,..., 90\%) and use root mean square pruning in the input layer to sparsify the layer to the target sparsity.  Sparsity is achieved by using a binary masking layer which can be represented mathematically the following way: 
\begin{equation}
    \mathcal{F}(\sss; \mathbf{W}) := R(\mathbf{W}^kR(\mathbf{W}^{k-1}\cdots \mathbf{M^1} \circ R(\mathbf{W^1 s}))),
\end{equation}
where $\mathbf{M^1}$ is a binary masking matrix that multiplies element-wise ($\circ$ represents the Hadamard product) to the input layer of the SDN.  The matrix $\mathbf{M^1}$ is inialized with ones.  During pruning, zeroing-out a node corresponds to making an entire column of $\mathbf{M^1}$ zeros.  During each prune, we zero-out the number of columns necessary to reach the target sparsity.  The pruned network is then retrained and the process repeats, incrementally increasing the target sparsity until very few inputs remain. To avoid overfitting to the training data, the remaining weights in the network are reinitialized after each prune.  The pruned inputs nodes have no effect on the prediction of the model because they remained zeroed-out throughout the entire process.  This is achieved by making the masking matrix non-trainable, meaning it cannot be updated during backpropagation.  \\
\indent Since the input nodes represent sensor measurements, we consider magnitude-based pruning of input nodes as a nonlinear method for sensor selection.  In contrast with the linear placement method presented in Section \ref{section:placement}, this algorithm, P-SDN, relies exclusively on SDNs to learn measurement selection and reconstruction in parallel.  


\section{Computational Experiments}
\label{section-experiments}
\indent Having introduced both Q-SDN and P-SDN, we compare the reconstructive performance of these algorithms to SDNs with randomly placed sensors and the more traditional modal reconstruction presented in Section \ref{section:placement} \cite{manohar-sensors}.  Performance is compared for two datasets, sea surface temperature \cite{sst-data} and the canonical example of fluid flow behind a cylinder at Reynolds number 100.  In both cases, reconstruction error (RE) is measured by
\begin{equation}
    \texttt{RE} = \frac{1}{N}\sum _{k = 1} ^{N} \frac{||\hat{\xx} _k - \xx_k ||}{||\xx_k||}
    \label{eq:vortexerror}
\end{equation}
where $\xx$ is the mean subtracted full state, $\hat{\xx}$ is the mean subtracted reconstructed state, and $k$ indexes the samples withheld during training.  The mean subtracted reconstructions are considered because in datasets where the mean state accounts for much of the spatial structure, the relative deviation provides a more sensitive measure of reconstructive performance \cite{erichson-decoder}.  We emphasize the importance of performance in the ``low sensor'' limit i.e. when increasing the number of sensors significantly improves reconstructive performance.  With larger numbers of sensors, the question of sensor placement becomes less important because even randomly chosen measurements can yield accurate reconstructions.  We find that Q-SDN consistently outperforms all other considered techniques.  Surprisingly, in Subsection \ref{subsection:p-sdn-sst}, we show P-SDN is significantly outperformed by SDNs with randomly chosen measurements (R-SDN).  R-SDN serves as an important benchmark for the performance of both Q-SDN and P-SDN as a successful placement scheme should outperform random placement.
\subsection{P-SDN and Sea Surface Temperature}
\label{subsection:p-sdn-sst}

We begin by considering the performance of P-SDN in comparison to that of R-SDN in an example dataset consisting of 1,400 samples of a a $180 \times 360$ grid representing sea surface temperature (44,219 points represent sea surface temperature, the rest landmass).  We initialize the P-SDN network with an input layer of 2,000 randomly chosen measurements, a first hidden layer of 350 nodes, a second hidden layer of 400 nodes, and an output layer of size 44,219.  The ADAM optimizer with learning rate of 0.001 is used to train all networks and an early stopping criteria with a patience of 5 is used. This early stopping criteria determines the intervals at which the network is pruned in our implementation of P-SDN. $20\%$ of the remaining nodes are pruned at each iteration and the resulting networks performances are compared to R-SDN with comparable numbers of sensors. In all cases, 1,000 samples are randomly selected to serve as training data, while the 400 withheld samples are used to evaluate performance.   We consider 32 instances of R-SDN and P-SDN to determine a mean and standard deviation of reconstruction error, as defined in Eq.  (\ref{eq:vortexerror}).  The results are displayed in Fig.  \ref{fig:p-sdn-fig}.

Surprisingly, R-SDN outperforms P-SDN at all numbers of sensors we consider; rather than learning ideal sensor locations, P-SDN as implemented here identifies poor locations.  Fig.  \ref{fig:pruning} illustrates why this might be the case.  While P-SDN identifies locations similar to that of QR, it restricts its placement to only a small subset of physical locations.  This appears to indicate that magnitude based pruning alone is not successful at identifying ideal placement, although we conjecture other pruning protocols might have greater success.  As a result of this finding, our subsequent discussion will focus on a comparison of Q-SDN to competing techniques.
\subsection{Q-SDN and Sea Surface Temperature}
\label{subsection:Q-SDN}

We now compare the performance of Q-SDN to that of R-SDN and modal reconstructions performed with sensor locations determined by the extension of Q-DEIM presented in \cite{manohar-sensors}.  Performance is considered in the same sea surface temperature dataset from the previous section and Q-SDN is implemented with the same network architectures and optimizer as presented in \ref{subsection:p-sdn-sst}.  The performance of R-SDN is evaluated by averaging the performance of 100 independent trials with varying, random sensor placement.  100 trials were also performed for both Q-SDN and the modal reconstruction, however measurement location was held constant in each instance.  The results, along with example reconstructions, are presented in Fig.  \ref{fig:sst_panel}.

Q-SDN uniformly outperforms both competing techniques.  Furthermore, although R-SDN produces comparably accurate reconstructions with sufficiently many sensors, in the low sensor limit the difference in accuracy between Q-SDN and R-SDN is pronounced. Modal reconstructions with principled sensor selection outperform R-SDN when very few sensors are available, but fail to improve in reconstruction accuracy with the addition of more sensors. 

\subsection{Q-SDN and Vortex Shedding Behind a Cylinder}\label{subsection-shedding}
\begin{figure}[t]
    \centering
    \includegraphics[width=0.65\columnwidth]{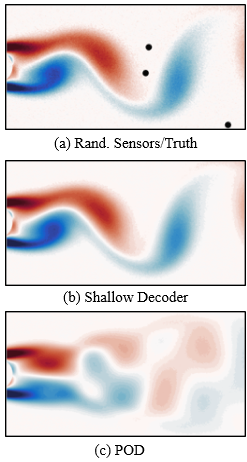}
    \caption{(a) Randomly placed sensors in the flow field of vortex shedding behind a cylinder, (b) reconstruction of flow field from the depicted sensors using a shallow decoder, and (c) reconstruction of flow field from the depicted sensors using standard POD methods.}
    \label{fig:randpod}
\end{figure}
\begin{figure*}[!t]
    \centering
    \includegraphics[width=\textwidth]{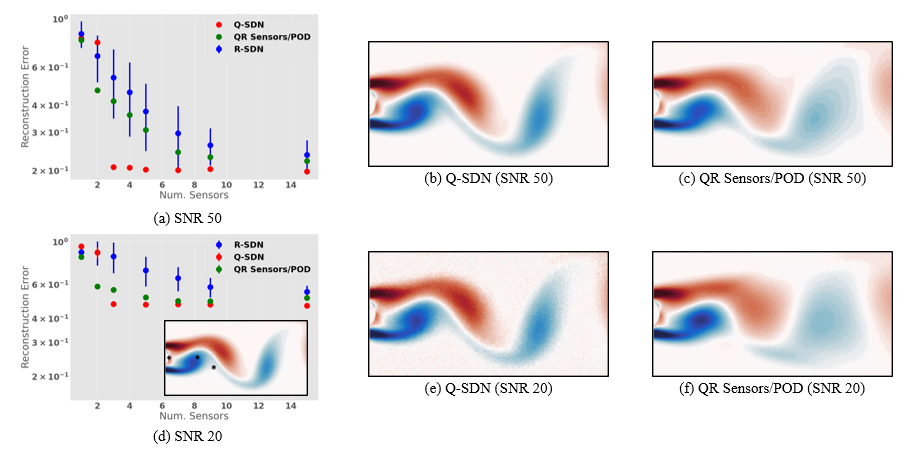}
    \caption{(a) Reconstruction errors obtained by Q-SDN, R-SDN, and POD combined with the placement from Section \ref{section:placement} at varying numbers of sensors with SNR 50, (b) reconstruction errors obtained for varying numbers of sensors at SNR 20 and QR sensor placement of 3 sensors within the true flow field, (b, c, e, f) reconstructions of the flow field depicted in (d).}
    \label{fig:vortex_shedding_panel}
\end{figure*}

The next system we consider is fluid flow behind a cylinder at Reynolds number 100. The training and test data used were synthetically generated by the immersed boundary projection method in two dimensions and consist of 100 and 51 samples, respectively.  Each complete state, $\xx$, is represented by a flattened $199 \times 384$ grid.  Since the simulation captures periodic behavior, the training data is highly representative of the test data.  We add Gaussian, random noise at varying peak signal-to-noise ratios (SNR) to simulate the corruption of the precise simulation.

As demonstrated by \cite{erichson-decoder}, SDNs are less sensitive to sensor placement than traditional, linear techniques, a result that is shown in Fig. \ref{fig:randpod}. With randomly placed sensors, R-SDN, the decoder network is able to reasonably reconstruct a particular sample in the withheld set while the POD based method fails to do so.  As a result, we consider the performance of Q-SDN in comparison to R-SDN and POD reconstructions with sensor locations determined by a QR decomposition and omit a comparison to POD reconstructions with random measurements. 

High-noise (SNR 20) and low-noise (SNR 50) regimes are considered with the number of available sensors ranging from 1 to 15. For both Q-SDN and P-SDN, the network architectures consist of an input layer with size equal to the number of available measurements, a first hidden layer with 35 nodes, a second hidden layer of 40 nodes, and an output layer of size 76,416, corresponding to each point in the $199 \times 384$ grid. Networks are trained with an ADAM optimizer with learning rate of 0.01 and an early stopping criteria with a patience of 5. POD based reconstructions are carried out as described in Section \ref{section:POD}, and the performance of all reconstructive techniques is evaluated by the metric presented in Eq. (\ref{eq:vortexerror}). The results of the computational experiments are presented in Fig. \ref{fig:vortex_shedding_panel}.

In both the high- and low-noise simulations, the performance of Q-SDN is superior than either R-SDN and the POD based reconstruction. However, the POD reconstructions typically outperform those obtained through R-SDN, in contrast to the sea surface temperature example. Q-SDN also appears to identify the optimal number of sensors to be 3, as in both the high- and low-noise instances there are marked improvements in reconstruction accuracy between the cases where 2 and 3 sensors are available, while the difference in performance between 3 and 15 sensors is negligible. 

\section{Discussion}
\label{section:discussion}

The importance of reconstructive tasks such as those considered here is ubiquitous across the physical sciences and engineering \cite{aerodynamics,manohar-sensors,ocean-modeling}. As a result, principled methodologies for both determining near optimal measurements and inferring complete state spaces from those measurements is of critical importance. We demonstrate that Q-SDN, a combination of a greedy linear selection with a nonlinear encoding, significantly outperforms competing reconstruction algorithms that either do not prescribe sensor placement (R-SDN), assign sensor placement by performing iterative magnitude based pruning (P-SDN), or rely on linear reconstructions (derivatives of gappy POD). In application to one example dataset (low-noise vortex shedding behind a cylinder), the reconstruction error of next best technique in the low sensor limit was nearly double that of Q-SDN.  Furthermore, to achieve comparably accurate reconstructions to that of Q-SDN with 3 sensors, R-SDN and linear reconstructions with optimized sensor locations required 5 times more sensors. Although the difference in performance was starkest when applied to the low-noise, vortex shedding dataset, Q-SDN also uniformly outperformed the considered competing techniques in application to the real-world dataset consisting of sea surface temperature data. We emphasize the importance performance with few sensors because of the numerous applications in which implementing many sensors is prohibitively costly or otherwise infeasible.  In such cases, Q-SDN shows promise to enable accurate reconstructions in regimes where other techniques fail.

The excellent performance of Q-SDN is made even more appealing by its relatively low computational cost and ease of implementation. QR decompositions, the basis of the placement mechanism used by Q-SDN, are natively implemented in most linear algebra packages and are an $\mathcal{O}(n^3)$ operation. The training protocol used here to learn a nonlinear mapping from a collection of measurements to full state spaces requires less than 5 minutes on a modern GPU. Although we only consider spatio-temporal systems here, Q-SDN can be applied more generally to reconstruct any image or set of measurements from some subsampling provided there is sufficient training data. As a result, we advocate the use of Q-SDN for principled reconstructions throughout the physical sciences. 

In application to the sea surface temperature data, the performance of the nonlinear sensor placement method, P-SDN, consistently performs worse than when sensors are randomly placed. The placement of the sensors may give some indication why this performance is so low. Sensor placement via magnitude-based pruning consistently clusters the sensors in cold regions of the ocean (See Fig. \ref{fig:qr_placement}). The temperature data, initially in degrees Celsius, is mean-subtracted and therefore the blue areas of the map indicate relatively large negative values. Magnitude-based pruning assumes that inputs connected to the first hidden layer by weak (i.e. small magnitude) weights are unimportant to the output of the network. After a cycle of pruning, the remaining inputs are the ones associated with the strongest weights. The weights connecting the inputs that represent cold regions of the ocean therefore have a larger magnitude than those connecting to warm temperature inputs. This indicates that simple magnitude-based pruning may not be suitable to this application, since it favors only the most negative inputs.   Ultimately, this reflects the need for further exploration for physics-informed pruning paradigms that may be able to make the pruning process more amenable to a given problem and application.

\section*{Acknowledgements}

The authors are especially indebted to Benjamin Erichson, Lionel Mathelin, Steven Brunton and Michael Mahoney for discussions concerning decoder networks.
JNK acknowledges support from the Air Force Office of Scientific Research (FA9550-19-1-0386 and FA9550-19-1-0011).

\bibliographystyle{ieeetr}
\bibliography{refs.bib}
\end{document}